\documentclass{amsart}

\usepackage{amsmath, amssymb, amsthm}
\usepackage{graphicx}
\usepackage{mathtools}
\usepackage{seqsplit}
\usepackage[english]{babel}
\usepackage[utf8]{inputenc}
\usepackage{float}
\usepackage{longtable}
\usepackage{marvosym}
\usepackage{xcolor}
\usepackage{ulem}
\usepackage[colorlinks=false]{hyperref}
\usepackage{cleveref}

\newtheorem{theorem}{Theorem}[section]
\newtheorem{lemma}[theorem]{Lemma}
\newtheorem{proposition}[theorem]{Proposition}

\theoremstyle{definition}
\newtheorem{definition}[theorem]{Definition}

\theoremstyle{remark}

\newtheorem{example}[theorem]{Example}

\definecolor{comment}{RGB}{102,0,102}

\crefname{theorem}{Theorem}{Theorems}
\crefname{lemma}{Lemma}{Lemmas}
\crefname{proposition}{Proposition}{Propositions}
\crefname{corollary}{Corollary}{Corollaries}
\crefname{definition}{Definition}{Definitions}
\crefname{example}{Example}{Examples}

\newcommand{\im}{\operatorname{Im}}

\begin{document}
\author{Jakub Szymański}
\address{Doctoral School of Exact and Natural Sciences, University of Warsaw, Poland}
\email{jt.szymanski@uw.edu.pl}
\thanks{The author was supported by the National Science Center Grant Maestro-13 UMO- 2021/42/A/ST1/00306.} % opcjonalnie

\title{Expansion of Integer Matrices over Various Rings}

% Kody MSC i słowa kluczowe
 % <-- wstaw swoje
          % <-- wstaw swoje

\begin{abstract} \sloppy
In this article, we explore the problem of constructing high-dimensional expanders through the study of relations between expansion constants over different rings. We investigate expansion constants of integer matrices regarded as morphisms between free modules over $\mathbb{R}$, $\mathbb{Z}$, and $\mathbb{Z}/p\mathbb{Z}$. We introduce a new condition which we call integral spanning regarding kernels of integer matrices, and prove that it ensures equality of real and integral expansions. In addition, we prove a bound on expansion constants over finite fields for a certain class of matrices in terms of the corresponding integral expansions. As an application, one may use this theorem to bound the expansion of codifferentials over $\mathbb{Z}/2\mathbb{Z}$ in degrees $0$ and $1$.
\end{abstract}

\maketitle

\section{Introduction}

Expander graphs are characterised by being sparse yet highly connected. Due to these two desirable traits, they have found applications in many areas of theoretical computer science. For example, see \cite{SipserSpielman1996} for their use in error-correcting codes, and \cite{AjtaiKomlosSzemeredi1983} for their use in sorting networks.

Given the central role of expanders in many constructions, it is natural to ask what higher-dimensional analogues might look like. It turns out that this generalisation can be introduced in many non-equivalent ways. There are definitions of spectral, coboundary, and cosystolic expanders. Nevertheless, all of them are in a similar vein and involve the augmented cellular cochain complexes of CW complexes $\left\{ X_n \right\}_{n \in \mathbb{N}}$, with coefficients in an integral domain $R$

\[
R \xrightarrow{} C^0_{\text{cell}}(X_n, R) \xrightarrow{} C^1_{\text{cell}}(X_n, R) \xrightarrow{} C^2_{\text{cell}}(X_n, R) \to \ldots 
\] 

The definition of high dimensional expansion then requires a uniform upper bound on expansion constants, as defined below, for all codifferential maps associated with the family $\left\{ X_n \right\}_{n \in \mathbb{N}}$, where, in the case of spectral expanders, we have $R = \mathbb{R}$ and for coboundary and cosystolic expansions $R = \mathbb{Z}/2\mathbb{Z}$.

\begin{definition}

    For a given morphism $A : R^n \longrightarrow R^m$, we define the expansion constant at a nonzero element $v \in \im(A)$ as

    \begin{align*}
        \Xi_R(A, v) 
        &= \inf \left\{ \frac{\| u \|}{\|v\|} \mid u \in A^{-1}(\{v\}) \right\},
    \end{align*}
    and a general expansion constant as

    \begin{align*}
        \Xi_R(A) &= \sup \left\{ \Xi_R(A,v) \mid v \neq 0 \in \im(A) \right\}.
    \end{align*}

\end{definition}

It turns out that the construction of spectral expanders is much simpler, and can be done by taking a family of finite covers of the presentation complex corresponding to a group with Property (T) (see \Cref{thm:T_and_expansion} \cite[Theorem 1.8]{uniform_waist}). This work is the result of considerations on cosystolic expansion constants, i.e., after change of field, in the case of the above sequence of finite covers. Namely, the principal question is about the relation of the expansion constants of integer matrices, considered as morphisms between free modules over $\mathbb{R}, \mathbb{Z}$ and $\mathbb{Z}/p\mathbb{Z}$ for a prime number~$p$. 

In \cite{uniform_waist} the authors prove that the equality between expansion over $\mathbb{Z}$ and $\mathbb{R}$ holds for totally unimodular matrices, i.e. such that all minors are either $-1,0$ or $1$. The authors also show that expansion behaves well with respect to exact sequences, provided the matrices are totally unimodular.

\begin{theorem}\cite[Theorem 2.8]{uniform_waist}\label{exact_unimod}
    Let $A$ and $B$ be integer matrices of size $m \times n$ and $k \times m$, such that the sequence $ \mathbb{R}^n \xrightarrow{A} \mathbb{R}^m \xrightarrow{B} \mathbb{R}^k $ is exact. If $A$ is totally unimodular, then $\Xi_{\mathbb{Z}}(A) = \Xi_{\mathbb{R}}(A)$ and $\Xi_{\mathbb{Z}}(B) = \Xi_{\mathbb{R}}(B)$.
\end{theorem}

In the context of codifferentials, the zeroth map is always totally unimodular, and the equality $\Xi_\mathbb{R}(\text{d}^0) = \Xi_\mathbb{Z}(\text{d}^0)$ holds for every CW complex. Moreover, if we add the condition about the trivial first cohomology group, based on the above theorem, we get the same equality for the first codifferential.

In the subsequent chapter we present a new condition which ensures the equality of real and integer expansions. Rather than insisting on the total unimodularity of the matrix, we establish it by formulating a condition on the kernel of the morphism defined by the matrix. We introduce this property for a free $\mathbb{Z}$-module in \Cref{int_completeness} and refer to the module as integrally spanned. Using more elementary methods, we prove an analogous statement of equal expansions over $\mathbb{R}$ and $\mathbb{Z}$ for the zeroth and first codifferentials, based on this condition.

Furthermore, we present a method for bounding the expansion over a finite field by the expansion constant over the integers. \Cref{finite_field_expansion} can be applied, for example, for the expansion of zeroth and first codifferentials over $\mathbb{Z}/2\mathbb{Z}$. Whether this result can be extended to a larger class of matrices remains an open problem. Addressing this question would make a significant step towards providing new examples of coboundary and cosystolic higher-dimensional expanders.

\section{New condition on equality of real and integer expansion}

In this chapter, we consider coefficient rings $R$ equal to $\mathbb{R}, \mathbb{Q}$, or $\mathbb{Z}$. For each of them, the norms on the free modules are defined as the sum of the absolute values of the coordinates of a vector. 

Our new condition is based on comparing finitely generated modules over integers with the corresponding modules over rationals generated by the same basis. The equalities we consider are required to hold for the generating sets after restriction to subsets of coordinates. The following definition introduces restriction maps \(p_I\), which can be viewed as projections onto the selected coordinates, preserving their order.

\begin{definition}
    Let $I = \{i_1 < i_2 <  \cdots  < i_m \}$ be a subset of $ \{1, \cdots ,n\}$ and let $R$ be a ring with unity. Then the morphism of $R$-modules $p_I : R^n \rightarrow R^m$ is defined by the formula
    \[
    p_I(x_1, \cdots ,x_n) = (x_{i_1}, \cdots ,x_{i_m}).
    \]
\end{definition}

\begin{definition}
    Let $R$ be a ring with unity and let $M$ be an $R$-module together with elements $m_1, \cdots, m_k \in M$. A module generated by them will be denoted by
    \[
    \langle m_1, \cdots, m_k \rangle_R \subseteq M.
    \]
\end{definition}

In the following, we shall not distinguish between \(p_I\) as morphisms of free modules over different rings \(R\), and we will use the same symbol \(p_I\) in all contexts.

\begin{definition}\label{int_completeness}
    Let $v_1, \cdots, v_k \in \mathbb{Q}^n$ be vectors with integer coefficients and let $V$ be a $\mathbb{Z}$-module generated by them, $V = \langle v_1, \cdots ,v_k \rangle_\mathbb{Z} \subseteq \mathbb{Q}^n$.
    We say that $V$ is integrally spanned if for each $I \subseteq \{1, \cdots,n\}$ the following equality of $\mathbb{Z}$-modules is satisfied
    \[
    \langle p_I(v_1), \cdots,  p_I(v_k) \rangle_\mathbb{Z} = \mathbb{Z}^{|I|} \cap \langle p_I(v_1), \cdots,  p_I(v_k) \rangle_\mathbb{Q}.
    \]
\end{definition}

% \textit{The independence of this property from the choice of generating set is routinely checked:} \\

\begin{lemma}
    The above property of a $\mathbb{Z}$-module does not depend on the choice of the generating set. 
\end{lemma}

\begin{proof}
Let us choose two generating sets $V = \{v_1, \cdots ,v_k\} \subseteq \mathbb{Z}^n \cap \mathbb{Q}^n$ and $W = \{w_1, \cdots, w_m\} \subseteq \mathbb{Z}^n \cap \mathbb{Q}^n$ such that the equalities from \Cref{int_completeness} are satisfied for $V$. We will show that they hold for $W$ as well. Since $V$ and $W$ are generating sets, we can introduce integers $\{a^i_{j}\}, \{b^i_{j}\}$ such that
\[
    v_i = \sum_{j = 1, \cdots ,m} a^i_j w_j, \quad w_i = \sum_{j = 1, \cdots , k} b^i_j v_j.
\]

Therefore, for rational numbers $\gamma_1, \cdots, \gamma_m \in \mathbb{Q}$ and a subset $I \in \{1, \cdots, n\}$ such that

\[
    w \coloneq \gamma_1 p_I(w_1) +  \cdots  + \gamma_m p_I(w_m) \in \mathbb{Z}^{|I|},
\]
the following holds

\[
    w = \left( \sum_{i = 1, \cdots ,m} \gamma_i b_1^i \right) p_I(v_1) +  \cdots  + \left( \sum_{i = 1, \cdots ,m} \gamma_i b_k^i \right) p_I(v_k),
\]
and thus $w =c_1p_I(v_1) +  \cdots  + c_kp_I(v_k)$ for some  $c_1, \cdots, c_k \in \mathbb{Z}$. Then we get

\[
    w = \left( \sum_{i = 1, \cdots ,k} c_ia^i_1 \right) p_I(w_1) +  \cdots  + \left( \sum_{i = 1, \cdots ,k} c_ia^i_m \right) p_I(w_m) \in \langle p_I(w_1), \cdots,  p_I(w_m) \rangle_\mathbb{Z},
\]
which proves that the equality

\[ 
    \langle p_I(w_1), \cdots,  p_I(w_m) \rangle_\mathbb{Z} = \mathbb{Z}^{|I|} \cap \langle p_I(w_1), \cdots,  p_I(w_m) \rangle_\mathbb{Q}
\]
holds for each $I \in \{1,\cdots,n\}$, which is exactly the condition from \Cref{int_completeness}.

\end{proof}

\begin{example}\label{example}
\leavevmode
\begin{enumerate}

    \item For a single vector $v \in \mathbb{Z}^n \cap \mathbb{Q}^n$ a  module $\langle v \rangle_\mathbb{Z}$ is integrally spanned if and only if all entries of $v$ are in $\{-1,0,1\}$. Otherwise, after restricting to an entry with a value $k \notin \{-1,0,1\}$, the equality from \Cref{int_completeness} does not hold
    \[
    \langle k \rangle_\mathbb{Z} \neq \mathbb{Z} = \mathbb{Z} \cap \langle k \rangle_\mathbb{Q}.
    \]
    
    \item For more than one vector, the set of possible entries is not restricted. Let $k \in \mathbb{Z}$ be any integer number, then the  $\mathbb{Z}$-module generated by the vectors $v_1 = (0,1), v_2 = (1,k)$ is integrally spanned.
    
    \item Extending the first example, a $\mathbb{Z}$-module generated by vectors $v_1, \cdots, v_k \in \mathbb{Z}^n \cap \mathbb{Q}^n$ is never integrally spanned if there exists a coordinate such that the corresponding entries of $v_1, \cdots, v_k$ have GCD greater than $1$. However, these conditions are not sufficient, as the example $v_1 = (1,1), v_2 = (1,3)$ shows, in which case
    \[
    \frac{1}{2}(1,1) + \frac{1}{2}(1,3) = (1,2) \notin \langle v_1,v_2 \rangle_\mathbb{Z},
    \]
    and of course it lies inside $\mathbb{Z}^2 \cap \langle v_1,v_2 \rangle_\mathbb{Q}$.

    \item Let $v_1, \cdots, v_k \in \mathbb{Q}^n$ be vectors with entries in $\{-1,0,1\}$ such that, for each coordinate $i \in \{1, \cdots ,n\}$ at most one vector has non-zero $i$-th entry. Then, for any rational numbers $\beta_1, \cdots, \beta_k$ and a subset $I \subseteq \{1, \cdots, n\}$, the linear combination $\beta_1 p_I(v_k) +  \cdots  +\beta_k p_I(v_k)$ has entries in $\{\pm\beta_1,  \cdots,  \pm \beta_k \}$. Therefore if such a combination belongs also to $\mathbb{Z}^{|I|}$, we may also represent it as a combination with integer coefficients. Thus, $\langle v_1, \cdots, v_k \rangle_\mathbb{Z}$ is integrally spanned.
    
\end{enumerate}
\end{example}

In \Cref{equlity_of_expansion_for_integral_spanned_matrix} we prove an important fact about the global minimum point of a function, which is linear over a finite number of components of the domain. The following definition is consistent with \cite{max-min}.

\begin{definition}
   Let $\Gamma$ be a closed convex set in $\mathbb{Q}^k$. A function $f: \Gamma \to \mathbb{Q}$ is said to be piecewise linear if there exists a finite family of closed sets whose union equals $\Gamma$ and $f$ coincides with some degree-one polynomial with rational coefficients on each of them.
\end{definition}

\begin{theorem}\label{equlity_of_expansion_for_integral_spanned_matrix}
    Let $A$ be an $m \times n$ integer matrix which we identify with the linear map $A : \mathbb{Z}^n \longrightarrow \mathbb{Z}^m$. If $\ker(A)$ is integrally spanned, the equality $\Xi_\mathbb{R}(A) = \Xi_\mathbb{Z}(A)$ of real and integral expansion holds.
\end{theorem}

Since the equality $\Xi_\mathbb{R}(A) = \Xi_\mathbb{Q}(A)$ follows from the density of rational numbers, we will focus on proving the equality $\Xi_\mathbb{Q}(A) = \Xi_\mathbb{Z}(A)$. To prove this, we need a lemma about finding the minimum of a piecewise linear function

\begin{lemma}\label{min_of_pwlin_func}
    Let $f: \mathbb{Q}^k \longrightarrow \mathbb {Q}$ be a function of the following form
    
    \[
    f(x_1,  \cdots , x_k) = \sum_{i = 1, \cdots ,n} \left| a_i + \sum_{j = 1, \cdots ,k} \xi_j^i x_j \right|,
    \]
    defined by the constants $\{a_i\} , \{ \xi_j^i\} \subset \mathbb{Q}$. Moreover, let 
    
    \begin{align*}
        H_i = \left\{ (x_1, \cdots ,x_k) \in \mathbb{Q}^k \mid a_i + \sum_{j = 1, \cdots ,k} \xi_j^i x_j = 0 \right\}, \\
        \mathcal{H} = \left\{ \bigcap_{i \in I} H_i \mid I \subseteq \{1, \cdots ,n\} \right\} \setminus \left\{ \emptyset \right\}. 
    \end{align*}

    Then, for any minimal element with respect to the inclusion $H \in \mathcal{H}$, the function $f$ is constant on $H$. Moreover, $f$ attains the global minimum at some affine subspace of this form.
\end{lemma}

\begin{proof}[Proof of \Cref{min_of_pwlin_func}]
    \sloppy
    The function assigning $(x_1, \cdots ,x_k)$ a value $\left| a_i + \sum_{j = 1, \cdots ,k} \xi_j^i x_j \right|$ is equal to $0$ on the affine subspace $H_i$ and linear on subsets which are formed from $\mathbb{Q}^k$ by diving it by this hyperplane. Therefore $f$ is piecewise linear and attains its minimum at the boundary of the subsets of the domain on which it is linear, which equals $\bigcup_{i \in \{1, \cdots ,n\}} H_i$. 
    
    %Since the boundaries of such subsets consist of elements of $\mathcal{H}$, on which $f$ is again piecewise linear, the point where $f$ attain its minimum belongs to some minimal element of $\mathcal{H}$. \\
    
    If we consider an element of $\mathcal{H}$ which is not minimal with respect to the inclusion
    
    \[
    H = \bigcap_{i \in I} H_i \in \mathcal{H},
    \]
    together with $\emptyset \neq I' \subseteq \{1, \cdots ,n\}$ such that $I \cap I' = \emptyset$ and $\forall_{l \in I'} H_l \cap H \neq \emptyset$ the function $f$ restricted to $H$ is again piecewise linear and boundary of subsets, on which it is linear, is given by the set $\bigcup_{l \in I'} H_l \cap H$. Therefore the minimum of $f_{|H}$ is reached somewhere at a smaller element of $\mathcal{H}$, considered with respect to the inclusion. This proves that the minimum value of $f$ is reached at some minimal element of $\mathcal{H}$.

    From now on, let $H = \bigcap_{i \in I} H_i \in \mathcal{H}$ be a minimal element for some $I \subseteq \{1, \cdots ,n\}$. By definition, $H$ is an affine space and we may describe it as $H = (y_1, \cdots ,y_k) + U$ for $(y_1, \cdots ,y_k) \in \mathbb{Q}^k$ and linear subspace $U \subset \mathbb{Q}^k$. 
    
    Let $l \in \{1, \cdots ,n\}$ do not belong to $I$ and let us consider the following linear function $w_l: U \longrightarrow \mathbb{Q}$
    
    \[
    w_l:  (x_1, \cdots ,x_k) \mapsto \sum_{j = 1, \cdots ,k} \xi_j^l x_j.
    \]
    Since 
    \[
    H \cap \left\{ (x_1, \cdots ,x_k) \in \mathbb{Q}^k \mid a_l + \sum_{j = 1, \cdots ,k} \xi_j^l x_j = 0 \right\} = \emptyset,
    \]
    the value $-\left(a_l + \sum_{j = 1, \cdots ,k} \xi_j^l y_j \right)$ does not belong to $\operatorname{Im}(w_l)$. Therefore $\operatorname{Im}(w_l) = 0$, which proves that for any point $(x_1, \cdots ,x_k)$ inside $H$ the value of $\left| a_l + \sum_{j = 1, \cdots ,k} \xi_j^l x_j \right|$ is the same. This proves that $f$ is constant on minimal elements of $\mathcal{H}$ and at one of them it attains its minimum value.
\end{proof}

\begin{proof}[Proof of \Cref{equlity_of_expansion_for_integral_spanned_matrix}]
   In what follows, we will also denote by $A$ a map between linear spaces over the rationals which is formally given by $A \otimes_\mathbb{Z} \mathbb{Q}$. However, it is desired to distinguish the kernels of both maps, denoted by $\ker_\mathbb{Z}(A)$ and $\ker_\mathbb{Q}(A)$, which are connected by the equation 
   \[
   \ker_\mathbb{Q}(A) =\ker_\mathbb{Z}(A) \otimes_\mathbb{Z} \mathbb{Q},
   \]
   since $\mathbb{Q}$ is a flat $\mathbb{Z}$-module.
    
    For every vector $v \in A(\mathbb{Q}^n)$ there exists an integer number $\beta$ such that $\beta v \in A(\mathbb{Z}^n)$ and hence
    
    \[
    \Xi_{\mathbb{Q}}(A) = \sup \left\{ \Xi_\mathbb{Q}(A,v) \mid v \in A(\mathbb{Z}^n)\right\}.
    \]
    In this regard, let us consider $v \in A(\mathbb{Z}^n)$, together with an integer vector $u = (u_1. \cdots ,u_n) \in A^{-1}(v)$. From the assumptions there exists a generating set of $\ker_\mathbb{Z}(A)$, consisting of integer vectors $\xi_1, \cdots ,\xi_k$ such that

    \[ 
    \forall_{I \in \{1. \cdots ,n\}} \quad \langle p_I(\xi_1), \cdots,  p_I(\xi_k) \rangle_\mathbb{Z} = \mathbb{Z}^{|I|} \cap \langle p_I(\xi_1), \cdots,  p_I(\xi_k) \rangle_\mathbb{Q}.
    \]
    From the definition of the expansion constant

    \[
    \Xi_\mathbb{Q}(A,v) = \frac{1}{||v||} \inf_{w \in A^{-1}(v)}||w||,
    \]
    and for each rational vector $w \in A^{-1}(v)$ there exist rational numbers $x_1, \cdots, x_k \in \mathbb{Q}$ such that $w = u + x_1 \xi_1 +  \cdots  + x_k\xi_k$. Thus,
    
    \[
    \Xi_\mathbb{Q}(A,v) = \frac{1}{||v||} \min_{x_1, \cdots, x_k \in \mathbb{Q}} \sum_{i = 1, \cdots ,n} \left| u_i + \sum_{j = 1, \cdots ,k} \xi_j^i x_j \right|,
    \]
    where $\xi_j^i$ denotes $i$-th coordinate of the vector $\xi_j$. By \Cref{min_of_pwlin_func} the minimum is reached at some point $(z_1, \cdots ,z_k) \in \mathbb{Q}^k$ belonging to an affine subspace 
    
    \[
    H = \bigcap_{i \in I} \left\{ (x_1, \cdots ,x_k) \in \mathbb{Q}^k \mid u_i + \sum_{j = 1, \cdots ,k} \xi_j^i x_j = 0 \right\},
    \] 
    for some $I \subseteq \{ 1, \cdots ,n \}$, on which the sum from the expansion definition is constant. This means that $z_1 p_I(\xi_1) +  \cdots  + z_kp_I(\xi_k)\in \mathbb{Z}^{|I|}$, which from the assumption of integral spanning implies the existence of integer values $\widetilde{z_1}, \cdots ,\widetilde{z_k} \in \mathbb{Z}$ such that $(\widetilde{z_1}, \cdots ,\widetilde{z_n}) \in H$ and therefore 
    \[
    \Xi_\mathbb{Q}(A,v) = 
    % \frac{1}{||v||} \min_{x_1, \cdots, x_k \in \mathbb{Z}} \sum_{i = 1, \cdots ,n} \left| u_i + \sum_{j = 1, \cdots ,k} \xi_j^i x_j \right| = 
    \Xi_\mathbb{Z}(A,v) \quad \text{and} \quad \Xi_\mathbb{Q}(A) = \Xi_\mathbb{Z}(A).
    \]
    
\end{proof}

Next, we consider a class of matrices arising from adjacency matrices of finite graphs and show that their kernels and images are integrally spanned. This, in turn, provides a new proof of Theorem 2.12 from the work \cite{uniform_waist}.

\begin{proposition}\label{int_dens_of_-1_0_1_matrix}
    Let $A$ be an $m \times n$ integer matrix whose entries are in $\{ -1, 0 ,1\}$ such that each row has at most one entry $1$ and one entry $-1$. Then the kernel and the image of the associated linear map $A: \mathbb{Z}^n \longrightarrow \mathbb{Z}^m$ are integrally spanned.
\end{proposition}

\begin{proof}
    Each row of $A$ is of one of the following types
    \begin{enumerate}
        \item is a zero vector,
        \item has only one nonzero value on $i$-th position, equal to 1 or -1,
        \item has two nonzero entries, $i$-th and $j$-th, which equal 1 and -1.
    \end{enumerate}

Let us associate with $A$ an undirected graph $G$ with $n$ vertices, labeled by $\{1, \cdots ,n\}$ and edges described by its rows. If one is of the second type, described above, we connect vertex $i$ to itself, and if a row is of the type $3$ we connect vertices $i$ and $j$. 

If $G$ has $k$ connected components we can describe them by vectors $v_1, \cdots ,v_k \in \mathbb{Z}^n$, where
\[
\text{$i$-th entry of } v_j =   \begin{cases}
            1 \quad \text{if $i$ is in $j$-th connected component,} \\
            0 \quad \text{otherwise}.
        \end{cases}
\]

The vectors $v_{i_1}, \cdots, v_{i_l}$, corresponding to those connected components for which no vertex is self-connected, are the basis of $\ker(A)$. Furthermore, for each $j \in \{1, \cdots ,n\}$, at most one vector among $v_{i_1}, \cdots, v_{i_l}$ has a non-zero $j$-th coordinate, which in that case equals $1$. As demonstrated in \Cref{example}, we can conclude that the kernel of $A$ is integrally spanned.

Next, $\im(A)$ is spanned by the vectors $A(e_1), \cdots, A(e_n)$, where $e_1, \cdots ,e_n$ is the standard basis of $\mathbb{Z}^n$ and for any rational linear combination of them the following equality holds
\[
\gamma_1 A(e_1) +  \cdots  + \gamma_n A(e_n) = A(\gamma_1, \cdots ,\gamma_n).
\]

Thus, if this combination is an integer vector, all entries corresponding to the indices belonging to the same connected component of the graph $G$ must have the same fractional part and if a connected component contains a self-connected vertex, the fractional part must equal 0. Therefore, the vector of fractional parts $(\{\gamma_1\}, \cdots ,\{\gamma_n\})$ belongs to the kernel and the integer linear combination satisfies
\[
(\gamma_1 - \{\gamma_1\})A(e_1) +  \cdots  + (\gamma_n - \{\gamma_n\}) A(e_n) = A(\gamma_1, \cdots ,\gamma_n).
\]

\sloppy

The images of the considered spanning vectors under $p_I$ are equal to $A_I(e_1), \cdots ,A_I(e_n)$, where $A_I$ denotes the matrix derived from $A$ by restricting it only to rows indexed by elements of $I$. In this case $A_I$ also satisfies the assumptions of this theorem, and the entire argument can be repeated. Thus $\im(A)$ is also integrally spanned.

\end{proof}

\begin{theorem}\cite[Theorem 2.12]{uniform_waist}\label{zero_homology_equal_exapansion}
    Let $X$ be a finite CW complex for which $H^1(X,\mathbb{R}) = 0$. Then the cellular differentials satisfy 
    \[
    \Xi_\mathbb{R}(d^0) = \Xi_\mathbb{Z}(d^0) \quad \text{and} \quad \Xi_\mathbb{R}(d^1) = \Xi_\mathbb{Z}(d^1)
    \]
\end{theorem}

\begin{proof}
    The matrix $d^0$ satisfies the assumptions of \Cref{int_dens_of_-1_0_1_matrix}. Therefore both the kernel and the image of $d^0$ are integrally spanned, and due to the triviality of the first cohomology group, $\ker(d^1) = \im(d^0)$. This result combined with the \Cref{equlity_of_expansion_for_integral_spanned_matrix} proves the equalities from the statement.
\end{proof}

\section{Expansion with $\mathbb{Z}_q$ coefficients}\label{fin_field_expansion}

Now, we focus on the case of coefficients given by a finite field $\mathbb{Z}_q \coloneq \mathbb{Z} / q\mathbb{Z}$ for a prime number $q$. Let us equip $\mathbb{Z}_q$ with the trivial absolute value

\[
|x| = \begin{cases} 
0 & \text{if } x = 0, \\ 
1 & \text{otherwise}. 
\end{cases}
\]
Due to the multiplicativity the above is the only possible absolute value for the rings considered.

\begin{theorem}\label{finite_field_expansion}
    Let $q$ be a prime number and let $A$ be an $m \times n$ matrix whose entries are in $\{ -1, 0 ,1 \}$ such that each row has at most one entry $1$ and one entry $-1$. Then for the corresponding linear map $\widetilde{A}: \mathbb{Z}_q^n \longrightarrow \mathbb{Z}_q^m$ the inequality $(q-1)\Xi_{\mathbb{Z}}({A}) \geq \Xi_{\mathbb{Z}_q}(\widetilde{A})$ holds.
\end{theorem}

\begin{proof}
    The finite field $\mathbb{Z}_q$ can be obtained from $\mathbb{Z}$ as a quotient ring. Therefore we have a natural projection $p$, which assigns to an integer number its remainder of division by $q$. We denote by $s$ the section of this projection given by
    \[
    \forall_{i \in \{0, 1,  \cdots , q-1\}} \quad [i] \overset{s}{\longmapsto} i.
    \]

    By $p$ and $s$ we also denote the induced maps between free products of the previously mentioned rings

    \[
    \mathbb{Z}^N \xlongrightarrow{p} \mathbb{Z}_q^N \xlongrightarrow{s} \mathbb{Z}^N.
    \]

    Let $w$ be a nonzero vector belonging to the image of $\widetilde{A}$ and let us denote 
    \[
    \alpha \coloneq \Xi_{\mathbb{Z}_q}(\widetilde{A}, w),
    \]
    together with $u \in \widetilde{A}^{-1}(w): \|u\|_{\mathbb{Z}_q^n} = \alpha \|w\|_{\mathbb{Z}_q^n}$. Then consider any vector $v \in \mathbb{Z}^n$ such that $Av = As(u)$. The projection applied to this equality yields

    \[
    \widetilde{A} p(v) = \widetilde{A} p(s(u)) = \widetilde{A} u = w.
    \]
    Hence, 
    
    \[
    \|v\|_{\mathbb{Z}^n} \geq \|p(v)\|_{\mathbb{Z}_q^n} \geq \alpha \|w\|_{\mathbb{Z}_q^n} \geq \alpha \|Av\|_{\mathbb{Z}^n} / (q-1),
    \]
    where the last inequality holds because each entry in the vector $Av = As(u)$ has an absolute value not greater than $q-1$ and $w = p(Av)$.
    Thus we obtain

    \begin{align*}
        (q-1) \Xi_{\mathbb{Z}}({A}, As(u)) \geq \Xi_{\mathbb{Z}_q}(\widetilde{A}, w),
    \end{align*}
    and

    \begin{align*}
        (q-1) \Xi_{\mathbb{Z}}({A}) \geq \sup_{u \in \mathbb{Z}_q^n} (q-1) \Xi_{\mathbb{Z}}({A}, As(u)) \geq \sup_{u \in \mathbb{Z}_q^n} \Xi_{\mathbb{Z}_q}(\widetilde{A}, \widetilde{A}u) = \Xi_{\mathbb{Z}_q}(\widetilde{A}).
    \end{align*}

\end{proof}

\section{Examples of various CW complexes with equal expansion constants}

For a finitely presented group $G = \langle s_1,  \cdots  ,s_n \mid r_1 \cdots ,r_m \rangle$ the presentation complex has homotopy group given by $G$ and the related augmented cellular cochain complex with coefficients in $\mathbb{Z}$ is

\[
\mathbb{Z} \xrightarrow{\text{d}_G^{-1}} \mathbb{Z} \xrightarrow{\text{d}_G^{0}} \mathbb{Z}^n \xrightarrow{\text{d}_G^{1}} \mathbb{Z}^m \xrightarrow{\text{d}_G^{2}} 0,
\]
where the matrix $\text{d}_G^0$ has zero coefficients and each row of $\text{d}^1_G$ is described by the vector counting the number of occurrences of generators in the corresponding relation (possibly with negative values).

We recall the following result on groups with property (T), which will be applied to the presentation complex described above.

\begin{theorem}\cite[Theorem 1.8]{uniform_waist} \label{thm:T_and_expansion}
    Let $X$ be a connected CW complex with a Kazhdan fundamental group and a finite 2-skeleton. There exists a constant $C > 0$ such that for every finite cover $\widetilde{X} \longrightarrow X$ the expansion constants of corresponding codifferentials $\Xi_{\mathbb{R}}(d^0_{\widetilde{X}})$ and $\Xi_{\mathbb{R}}(d^1_{\widetilde{X}})$ are bounded by $C$.
\end{theorem}

In this section, two classes of groups are to be considered

\begin{enumerate}
    \item Let us consider be the braid groups on $n$ strands 
    \begin{align*}
           B_n = \langle \sigma_1, \sigma_2, \dots, \sigma_{n-1} \mid
           & [\sigma_i ,\sigma_j] \text{ for } |i-j| \geq 2, \ \\
           & \sigma_i \sigma_{i+1} \sigma_i \sigma_{i+1}^{-1} \sigma_i^{-1} \sigma_{i+1}^{-1} \text{ for } i < n-1\rangle.
    \end{align*}
    Based on this presentation, $d^1_{B_n}$ has rows given by zero vectors for the rows corresponding to the relations of the first type, or vectors with $1$ on entry $i$ and $-1$ on entry $i+1$ in the second case. Due to \Cref{equlity_of_expansion_for_integral_spanned_matrix} and \Cref{finite_field_expansion} the following holds
    \[
    \Xi_\mathbb{R}(\text{d}^1_{B_n}) = \Xi_\mathbb{Z}(\text{d}^1_{B_n}) \geq \Xi_{\mathbb{Z}_2}(\text{d}^1_{B_n}).
    \]

    \item Next we focus on the Steinberg groups
    \begin{align*}
    St_n(\mathbb{Z}) = \langle x_{ij} \text{ for } 1 \leq i \neq j \leq n \mid
     & [x_{ij}, x_{kl}] = 1 \text{ for } i \neq l, j \neq k, \\
    & [x_{ij}, x_{jk}] = x_{ik} \text{ for } i \neq k \rangle.
    \end{align*}
    Again, $d^1_{St_n(\mathbb{Z})}$ satisfies the conditions of proved theorems and the same relations between expansion constants holds as in the case of the braid groups. What is more, due to \cite{ershov2010property}, Steinber groups for $n \geq 3$ have property (T). This means that by \Cref{zero_homology_equal_exapansion} and \Cref{thm:T_and_expansion} for a fixed $n \geq 3$ and any finite cover $\widetilde{X} \longrightarrow X_{St_n(\mathbb{Z})}$ we have the following equalities

    \[
    \Xi_\mathbb{R}\left(d^0_{\widetilde{X}}\right) = \Xi_\mathbb{Z}\left(d_{\widetilde{X}}^0\right), \quad \Xi_\mathbb{R}\left(d_{\widetilde{X}}^1\right) = \Xi_\mathbb{Z}\left(d_{\widetilde{X}}^1\right),
    \]
    and these expansion constants are uniformly bounded for all finite covers of the presentation complex.

\end{enumerate}

%% References


\bibliographystyle{alpha}
\bibliography{references.bib}

\end{document}